\documentclass[11pt]{article}
\usepackage{latexsym}
	\usepackage{amscd}
	\usepackage{amssymb}
	\usepackage{amsmath}
	\usepackage{amsthm}
	\usepackage{enumerate}
	\usepackage{verbatim}
\def\co{\colon\thinspace}

\newtheorem{thm}{Theorem}
\newtheorem{prop}[thm]{Proposition}

\newtheorem{Example}[thm]{Example}
\newenvironment{ex}{\begin{Example}\rm}{\end{Example}}
\newtheorem{Counterexample}[thm]{Counterexample}

\newtheorem{remark}[thm]{Remark}

\newtheorem{Fact}[thm]{Fact}

\newtheorem{Nothing}[thm]{$\!\!\!$}

\newcommand{\be}{\begin{equation}}
\newcommand{\ee}{ \end{equation}} 
\newcommand{\ba}{\begin{eqnarray}}
\newcommand{\ea}{\end{eqnarray}}
\newcommand{\ban}{\begin{eqnarray*}}
\newcommand{\ean}{\end{eqnarray*}}

\begin{document}
\abovedisplayskip=6pt plus3pt minus3pt
\belowdisplayskip=6pt plus3pt minus3pt
\title{\bf Vector bundles with infinitely many souls
\rm}
\author{Igor Belegradek}
\date{}
\maketitle
\begin{abstract}\noindent
We construct the first examples of manifolds,
the simplest one being $S^3\times S^4\times\mathbb R^5$,
which admit infinitely many complete nonnegatively curved 
metrics with pairwise nonhomeomorphic souls.
\end{abstract}
%\tableofcontents

According to the soul
theorem of J.~Cheeger and D.~Gromoll~\cite{CG}, a complete open manifold of
nonnegative sectional curvature is diffeomorphic to the total space of the
normal bundle of a compact totally geodesic submanifold, called a soul.
The soul is not unique but any two souls are mapped to each other
by an ambient diffeomorphism inducing an isometry
on the souls~\cite{Sha}.
In this note we show that the homeomorphism type of the
soul generally depends on the metric, namely the following is true.

\begin{thm}\label{s3s4r5} There exist infinitely many complete Riemannian
metrics on $S^3\times S^4\times\mathbb R^5$  
with $\sec\ge 0$ and pairwise nonhomeomorphic souls.
\end{thm}

The proof applies some classical techniques of geometric 
topology to recent examples of nonnegatively curved
manifolds due to K.~Grove and W.~Ziller~\cite{GZ}.
As we explain below, it is much easier
to produce a manifold with {\it finitely many} nonnegatively curved metrics 
having nonhomeomorphic souls, 
however the full power of~\cite{GZ} is needed to get
infinitely many such metrics.

Grove and Ziller~\cite{GZ} showed
that any principal $S^3\times S^3$-bundle
over $S^4$ admits an $S^3\times S^3$-invariant metric
with $sec\ge 0$. By O'Neill's formula, all associated
bundles admits metrics with $\sec\ge 0$ which
gives rise to a rich class of examples, including
all sphere bundles over $S^4$ with structure group $SO(4)$.
Note that the souls in Theorem~\ref{s3s4r5} are the total spaces of 
$S^3$-bundles over $S^4$ with structure group $SO(3)$.
Theorem~\ref{s3s4r5} is a particular case of the following.

\begin{thm}\label{main thm}
Let $\xi$ be a rank $n$ 
vector bundle over $S^4$ with structure group $SO(3)$, 
let $q\co S\to S^4$ be a smooth $S^{m-1}$-bundle with
structure group $SO(3)$, 
let $\eta$ be the $q$-pullback of $\xi$. 
If $m=4$, $n>4$, or if $m>4$, $n>m+3$, then the total
space of $\eta$ admits infinitely many complete Riemannian
metrics with $\sec\ge 0$ and pairwise nonhomeomorphic souls.
\end{thm}

The main topological tool used in this paper is 
a result of L.~Siebenmann~\cite{Sie}
that generalizes the famous Masur's theorem: any tangential 
homotopy equivalence of closed smooth $n$-manifolds is 
homotopic to a diffeomorphism after taking the product with 
the identity map of $\mathbb R^{n+1}$. 
Exotic $7$-spheres are stably parallelizable, so they all
become diffeomorphic after taking the product with $\mathbb R^8$,
and in fact it suffices to take product with $\mathbb R^5$.
A homotopy $7$-sphere is called a {\it Milnor sphere} if
it is diffeomorphic to the total space of an $S^3$-bundle
over $S^4$; 
it is known that $10$ out of $14$ homotopy $7$-spheres are Milnor.
According to~\cite{GZ}, any Milnor sphere carries a metric 
with $\sec\ge 0$, which leads to the following.

\begin{prop}\label{exotic} For every Milnor sphere $\Sigma$, 
the manifold $S^7\times \mathbb R^5$ has a complete Riemannian
metric of $sec\ge 0$ with soul diffeomorphic to $\Sigma$.
\end{prop}

Gromoll and Tapp recently classified~\cite{GT} all nonnegatively curved 
metrics on $S^2\times\mathbb R^2$, and asked for a similar
classification on $S^n\times \mathbb R^k$.
The above corollary indicates that such a classification
would be rather involved.

It would be interesting to classify the total spaces
of $S^3$-bundles over $S^4$ with nonzero Euler class up to
tangential homotopy equivalence. Indeed, by Masur's theorem 
such a classification should lead to an analog of 
corollary~\ref{exotic}, with $S^7$ replaced by any 
$S^3$-bundle $S$ over $S^4$ with nonzero Euler class.
It is worth mentioning that non-vanishing of the Euler class 
implies that the homotopy type of $S$ contains at most finitely
many nondiffeomorphic $S^3$-bundles over $S^4$~\cite{Tam} 
(cf.~\cite{CE}). 

It is easy to construct manifolds admitting
two nonnegatively curved metrics with pairwise
nonhomeomorphic souls. Here we describe two such situations
where souls are lens spaces, or simply-connected
homogeneous manifolds.
\begin{ex}
%\bf (Souls are lens spaces)\rm\ \
It is well-known~\cite[pp. 96, 100]{Coh} that 
the lens spaces $L(7,1)$ and $L(7,2)$ are homotopy equivalent,
but not simply-homotopy equivalent (hence nonhomeomorphic).
Orientable $3$-manifolds are parallelizable,
so any homotopy equivalence of $L(7,1)$, $L(7,2)$ is tangential,
and therefore~\cite[Theorem 2.3]{Sie}, we get a diffeomorphism 
of manifolds $L(7,1)\times\mathbb R^4$ and $L(7,2)\times\mathbb R^4$
admitting obvious product metrics of $\sec\ge 0$ with souls
$L(7,1)\times\{0\}$, $L(7,2)\times\{0\}$.
Since every homotopy type contains at most finitely many
nonhomeomorphic lens spaces, this procedure yields only 
finitely many pairwise nonhomeomorphic souls 
on a given manifold.
\end{ex}

\begin{prop} \label{homog}
If $i$ is a positive integer is divisible by $24$, then
there is a compact homogeneous space $G/H$ homotopy equivalent
but not homeomorphic to $S^{4i-1}\times S^4$, such that
for $n>4i+3$, the manifold
$G/H\times\mathbb R^n$ carries two complete nonnegatively curved
metrics with souls diffeomorphic to $G/H$ and $S^{4i-1}\times S^4$.
\end{prop}

The above mentioned result of Siebenmann~\cite[Theorem 2.2]{Sie} 
states that if $\xi$ is a rank $>2$ vector bundle over a compact
manifold, and the total space $E(\xi)$ of $\xi$ contains
a smooth compact submanifold $S$ such that
the inclusion $S\to E(\xi)$ is a homotopy equivalence,
then $E(\xi)$ has a vector bundle structure with zero 
section $S$. 

Then it easily follows~\cite[Theorem 2.3]{Sie}
that if $t$ is tangential homotopy equivalence
of $n$-manifolds, each being the total space of a vector
bundle over a compact smooth manifold of dimension less than
the rank of the bundle, then $t$ is homotopic to
a diffeomorphism. 
The upshot is that the total space of vector bundle
of sufficiently large rank often has many other
vector bundle structures, perhaps with nonhomeomorphic 
base spaces, and this fact plays a crucial role in this 
note.

\begin{proof}[Proof of Theorem~\ref{main thm}]
In what follows we denote the total space of a vector bundle $\zeta$ 
by $E(\zeta)$, and the associated sphere bundle by $S(\zeta)$.
Principal $SO(3)$-bundles over $S^4$ are in one-to-one
correspondence with $\pi_4(BSO(3))\cong\mathbb Z$.
Let $P_k$ be the principal $SO(3)$-bundle over $S^4$
corresponding to $k\in\pi_4(BSO(3))$. 
Let $\xi_k^n$ be the rank $n$ vector bundle over $S^4$ associated
with $P_k$ via the standard inclusion $SO(3)\to SO(n)$.

Let $P_{k,l}$ be the principal $SO(3)\times SO(3)$-bundle
over $S^4$ which is the pullback of $P_k\times P_{l}$ via
the diagonal map $\triangle\co S^4\to S^4\times S^4$. 
According to~\cite{GZ}, 
$P_{k,l}$ admits an $SO(3)\times SO(3)$-invariant metric
with $sec\ge 0$. 
Consider the $S^{m-1}\times\mathbb R^n$-bundle 
\[\triangle^\#(S(\xi_k^m)\times\xi_{l}^n)=
S^{m-1}\times_{SO(3)\times 1} P_{k,l}\times_{1\times SO(3)}\mathbb R^n\] 
associated with $P_{k,l}$.
The total space of $\triangle^\#(S(\xi_k^m)\times\xi_{l}^n)$ 
is also the total space of a rank $n$ vector bundle 
over $S(\xi_k^m)$ which we denote by $\eta^{m,n}_{k,l}$.
Note that $\eta^{m,n}_{k,l}$ is the pullback of $\xi_{l}^n$ 
via the projection $q^m_k\co S(\xi_k^m)\to S^4$.
By O'Neill's formula for submersions, $E(\eta^{m,n}_{k,l})$ carries
a complete metric with $\sec\ge 0$ with zero section being a soul.

First, note that $S(\xi_k^m)$ is fiber homotopy equivalent to 
$S(\xi_i^m)$ if  $k\equiv i\ (\mathrm{mod}\ 12)$. Indeed,
$S^{m-1}$-fibrations over $S^4$ are classified, up to
fiber homotopy equivalence, by $\pi_3(SG_{m})$
where $SG_{m}$ is the space of orientation-preserving
self-homotopy equivalences of $S^{m-1}$.
The fibrations $S(\xi_k^m)$ are classified by the image
of \[\phi\co\pi_3(SO(3))\to \pi_3(SO(m))\to\pi_3(SG_{m}).\]
Since $m\ge 4$, $\phi$ factors as
$\pi_3(SO(3))\to \pi_3(SF_3)\to\pi_3(SG_{4})\to \pi_3(SG_{m})$ 
where $SF_3$ is the space
of base-point-preserving elements of $SG_4$.
It is well-known that $SF_3$ is the identity component
of the loop space $\Omega^3 S^3$~\cite[Chapter 3]{MM}, hence 
$\pi_3(SF_3)\cong\pi_6(S^3)$.
Thus, $\phi$ factors through the $J$-homomorphism
$\mathbb Z\cong\pi_3(SO(3))\to \pi_6(S^3)\cong\mathbb Z_{12}$,
and the result follows.  

Second, show that $S(\xi_k^m)$ is homeomorphic to $S(\xi_i^m)$ iff $k=\pm l$
iff $S(\xi_k^m)$ is diffeomorphic to $S(\xi_l^m)$.
Indeed, the tangent bundle of $S(\xi_k^m)$ is stably isomorphic to
the $q_k^m$-pullback of $\xi_k^m$, because $S(\xi_k^m)$ is a two-sided
hypersurface in $E(\xi_k^m)$. 
Since $\xi_k^m$ is an $SO(3)$-bundle,
the first Pontrjagin class of $\xi_k^m$ is the $\pm 4k$-multiple
of a generator of $H^4(S^4,\mathbb Z)$ (cf.~\cite{Mil}). 
Also $S(\xi_k^m)$
has a section so that $q_k^m$ induces an isomorphism 
on the $4$th cohomology, hence $p_1(TS(\xi_k^m))$
is the $\pm 4k$-multiple of a generator.
By the topological invariance of rational Pontrjagin classes,
$S(\xi_k^m)$ is not homeomorphic to $S(\xi_i^m)$ unless
$k=\pm i$. Finally, $P_k$ is the pullback of $P_{-k}$ via
an orientation-reversing self-diffeomorphism of $S^4$,
so $S(\xi_k^m)$ and $S(\xi_{-k}^m)$ are diffeomorphic. 

Third, note that $E(\eta^{m,n}_{k,l})$ and  $E(\eta^{m,n}_{i,j})$
are tangentially homotopy equivalent provided
$k\equiv i\ (\mathrm{mod}\ 12)$, and $k+l=i+j$.
Indeed, fix $k,l,i,j$ with these properties.
The tangent bundle to $E(\eta^{m,n}_{k,l})$ is 
determined by its restriction to the zero section.
This restriction is the
$q_k^m$-pullback of $\xi_k^m\oplus\xi_l^n$ for any $k,l$. 
The images of $\xi_k^m$, $\xi_{l}^n$
under the homomorphism 
\[\pi_4(BSO(3))\to\pi_4(BSO)\] add up
to $k+l$, and the addition in $\pi_4(BSO)$ is given by
the Whitney sum $\oplus$. Thus $k+l=i+j$ implies that 
$\xi_k^m\oplus\xi_l^n$ and $\xi_i^m\oplus\xi_j^n$
are stably isomorphic.
Since $k\equiv i\ (\mathrm{mod}\ 12)$ there is a fiber
homotopy equivalence $g\co S(\xi_k^m)\to S(\xi_i^m)$, 
so that $g\circ q_i^m$ is homotopic to $q_k^m$.
Therefore, $g$ induces a tangential homotopy equivalence 
$t\co E(\eta^{m,n}_{k,l})\to E(\eta^{m,n}_{i,j})$, which 
is the composition of the projection of $\eta^{m,n}_{k,l}$,
followed by $g$, and then by the zero section of 
$\eta^{m,n}_{i,j}$.

Next we show that $E(\eta^{m,n}_{k,l})$ and
$E(\eta^{m,n}_{i,j})$ are diffeomorphic
if $k\equiv i\ (\mathrm{mod}\ 12)$, and $k+l=i+j$.
By Haefliger's embedding theorem~\cite{Hae},
the restriction of $t$ to the zero section of 
$\eta^{m,n}_{k,l}$ is homotopic to a smooth embedding $f$
because $2n\ge m+6$, which means we are in metastable range.
Since $n\ge 3$,~\cite[Theorem 2.2]{Sie} implies that
$E(\eta^{m,n}_{i,j})$ has a vector bundle structure
with zero section $f$. 
Since $t$ is tangential, $\eta^{m,n}_{k,l}$ is stably
isomorphic to $\nu_f$ which is the normal bundle to $f$.

In fact, our assumptions on $n,m$ imply that 
$\eta^{m,n}_{k,l}$ and $\nu_f$ are isomorphic.
Indeed, if $m>4$, $n>m+3=dim(S(\xi_k^m))$, then
we are in stable range, hence $\eta^{m,n}_{k,l}\cong\nu_f$.
If $m=4$, $n>4$, we apply
obstruction theory comparing $\nu_f$, $\eta^{m,n}_{k,l}$,
which are thought of as classifying maps from $S(\xi_k^m)$ to 
$BSO(n)$. Since $S(\xi_k^m)$ has a section, it can be obtained
by attaching a $7$-cell to $S^{3}\vee S^4$.
The bundles  $\nu_f$, $\eta^{m,n}_{k,l}$
are isomorphic on the $6$-skeleton $S^{3}\vee S^4$, 
because they are stably isomorphic and $n>4=\dim (S^{3}\vee S^4)$.
Comparing $\nu_f$, $\eta^{m,n}_{k,l}$ on the $7$-cell, 
we get a map $S^7\to BSO(n)$ which is nullhomotopic
since $\pi_7(BSO(n))=0$ if $n>4$ \cite[page 970]{Mim}).
Thus, $\nu_f$, $\eta^{m,n}_{k,l}$ are isomorphic, and
hence  $E(\eta^{m,n}_{k,l})$ and
$E(\eta^{m,n}_{i,j})$ are diffeomorphic.

To summarize, each manifold
$E(\eta^{m,n}_{k,l})$ has infinitely many vector bundle
structures $\eta^{m,n}_{i,j}$
with base manifolds $S(\xi^m_i)$ for any $i,j$ satisfying 
$k\equiv i\ (\mathrm{mod}\ 12)$, and $k+l=i+j$,
and the proof is complete.
\end{proof}

\begin{proof}[Proof of Proposition~\ref{exotic}]
Let $\Sigma$ be a homotopy $7$-sphere of $\sec\ge 0$~\cite{GM, GZ}.
The product metric on $\Sigma\times\mathbb R^n$ has $sec\ge 0$
with a soul $\Sigma\times\{0\}$.
By~\cite{Hae} any homotopy equivalence $S^7\to \Sigma\times\mathbb R^n$ 
is homotopic to a smooth embedding if $n\ge 5$,
so $\Sigma\times\mathbb R^n$ gets a structure of a vector
bundle over $S^7$ which is necessarily trivial
since $\pi_7(BSO(n))=0$ for $n\ge 5$.
Thus, $\Sigma\times\mathbb R^5$ is diffeomorphic to
$S^7\times\mathbb R^5$, as promised. 
\end{proof}

\begin{proof}[Proof of Proposition~\ref{homog}]
Kamerich~\cite[page 116]{Kam} (cf.~\cite[page 275]{Oni}) 
showed that if $i$ is a positive integer divisible by $24$, then
there is a compact homogeneous space $G/H$ homotopy equivalent
but not homeomorphic to $S^{4i-1}\times S^4$. 
Here $G=Sp(i)\times Sp(2)$ and $H=Sp(i-1)\times Sp(1)\times Sp(1)$,
where $Sp(i-1)$ is embedded in $Sp(i)$ is the standard way,
the first $Sp(1)$ is embedded into $Sp(i)\times Sp(2)$ diagonally 
so that a quaternion goes to the last diagonal entry of the
matrix in $Sp(i)$, and also to the first diagonal 
entry in $Sp(2)$, while the second $Sp(1)$ 
goes into the last diagonal entry of $Sp(2)$.

Since $n>4i+3$, 
the homotopy equivalence $S^{4i-1}\times S^4\to G/H\times\mathbb R^n$
is homotopic to a smooth embedding $f$. 
By~\cite[Theorem 2.2]{Sie}, $G/H\times\mathbb R^n$
admits an $\mathbb R^n$-bundle structure with zero section $f$.

Note that any $\mathbb R^n$-bundle $\xi$ over $S^{4i-1}\times S^4$
is the pullback of a $\mathbb R^n$-bundle over $S^4$ via the projection
$S^{4i-1}\times S^4\to S^4$. This is proved by obstruction theory
for maps $S^{4i-1}\times S^4\to BSO$ by comparing $\xi$ with 
the bundle $\xi_4$ obtained by pullbacking $\xi$ to $S^4$
via an inclusion $S^4\to S^{4i-1}\times S^4$, and then pullbacking it
back to  $S^{4i-1}\times S^4$ via the projection  
$S^{4i-1}\times S^4\to S^4$. Since $\pi_{4i-1}(BSO)=0$, 
$\xi$ and $\xi_4$ agree on  $S^{4i-1}\vee S^4$, and they are 
homotopic on the top $4i+3$-cell as $\pi_{4i+3}(BSO)=0$.

Since each vector bundle over $S^4$ 
carries $\sec\ge 0$ with zero section being a soul~\cite{GZ},
so does the product of the bundle and $S^{4i-1}$. Thus
$G/H\times\mathbb R^n$ gets a metric with $\sec\ge 0$ and soul
$S^{4i-1}\times S^4$. On the other hand, $G/H\times\mathbb R^n$ 
has the product metric with soul $G/H$.
\end{proof}

\bf Acknowledgements.\rm\ The results of this paper
were obtained during a geometry meeting in 
Oberwolfach in June 2001, and the author is grateful to
the meeting organizers for excellent working conditions
and hospitality. Special thanks are due to Wolfgang Ziller
for insightful discussions and encouragement, and
to McKenzie Wang for helpful communications
about Kamerich's thesis.

\small
\bibliographystyle{amsalpha}
\bibliography{soul-new}

\

DEPARTMENT OF MATHEMATICS, 253-37, CALIFORNIA INSTITUTE OF TECHNOLOGY,
PASADENA, CA 91125, USA

{\normalsize
{\it email:} \texttt{ibeleg@its.caltech.edu}}

\end{document}